\documentclass[12pt]{amsart}
\usepackage[utf8]{inputenc}
\usepackage{amssymb}
\usepackage{amsmath}
\usepackage{amsthm}
\usepackage{graphicx}
\usepackage{amscd}
\usepackage{epic}
\usepackage[alphabetic]{amsrefs}
\usepackage{xypic}
\usepackage[colorlinks, citecolor=magenta]{hyperref}
\usepackage[left=1.5in,right=1.5in,top=1.5in,bottom=1.5in]{geometry}
\usepackage{soul}
\usepackage{tikz-cd}
\usepackage[pagewise]{lineno}%\linenumbers
\usepackage{comment}
\usepackage{enumitem}
\usepackage{mathrsfs}
\usepackage{tikz-cd}
\theoremstyle{plain}
\newtheorem{theorem}{Theorem}[section]

\newtheorem{lemma}[theorem]{Lemma}

\newtheorem{definition-lemma}[theorem]{Definition-Lemma}
\theoremstyle{remark}
\newtheorem{remark}[theorem]{Remark}

\theoremstyle{definition}

\newcommand{\supp}[0]{\operatorname{Supp}}

\def\>{\geq}
\newcommand{\mbQ}{\mathbb{Q}}
\newcommand{\mbR}{\mathbb{R}}

\def\mcO{\mathcal{O}}

\newcommand{\mult}{{\rm mult}}
\newcommand{\Supp}{{\rm Supp}}

\newcommand{\mbC}{\mathbb{C}}

\def\mbP{\mathbb{P}}
\def\>{\geq}

\def\mcO{\mathcal{O}}

\def\Nklt{\operatorname{Nklt}}
\def\Nlc{\operatorname{Nlc}}
\def\mult{\operatorname{mult}}

\def\dim{\operatorname{dim}}

\theoremstyle{definition}
\newtheorem{definition}[theorem]{Definition}
\theoremstyle{definition}

\numberwithin{equation}{section}

\theoremstyle{remark}
\newtheorem{claim}[theorem]{Claim}

\author{Swapnajit Das}
\address{School of Maths, Tata Institute of Fundamental Research, Mumbai}
\email{das@math.tifr.res.in}

\title[Connectedness Principle]{Connectedness Principle in the Complex Analytic Setting}
\begin{document}
\begin{abstract}
We show that if $(X,B)$ is a normal compact K\"ahler pair and $f:X\to Z$ is a proper morphism with connected fibers to a normal compact K\"ahler variety $Z$ such that $(-K_X+B)$ is nef over $Z$, then $\Nklt(X,B)\cap f^{-1}(z)$ has at most two connected components for every $z\in Z$. More generally, we prove the result for generalized pairs whose nef part is allowed to be transcendental.  Moreover, we record that the same result holds for projective morphism between analytic varieties (possibly non-compact).
\end{abstract}

\maketitle
\tableofcontents

\section{Introduction}
Pairs and their singularities play a fundamental role in higher dimensional algebraic geometry. The connectedness of the non-klt locus of a log pair is an important topic in higher dimensional geometry. Let us consider a very simple example: let $X=\mbP^1$ and $B=p_1+p_2+\cdots p_k$ where $p_1,p_2,\dots,p_k$ are distinct points. Note that $\Nklt(X,B)$ can have any number of connected components. Similarly, in higher dimensions, one can consider a smooth projective variety with a simple normal crossing boundary divisor and show that the non-klt locus can have any number of connected components. However, in special situations, the non-klt locus shows interesting behavior. For example, Shokurov proved that if $X$ is a surface and $-(K_X+B)$ is ample, then $\Nklt(X,B)$ is connected (see \cite[Lemma 5.7]{Sho93}). Later, Koll\'ar generalized this result to higher dimensions (\cite[Theorem 17.4]{Kol92}). In fact, the result holds if we only assume $-(K_X+B)$ is nef and big. It also holds in the relative setting: if $X\to S$. is a morphism between projective varieties with connected fibers, the non-klt locus is connected near every fiber.\\

A natural question is what happens if we only assume $-(K_X+B)$ is nef. A simple example $(X=\mbP^1, B=p_1+p_2)$ shows that $K_X+B\equiv 0$ but $\Nklt(X,B)$ has two connected components. Hacon and Han (see \cite{HH19}) investigated the above question. They proved that if $\dim X \leq 4$ and $-(K_X+B)$ is nef, then $\Nklt(X,B)$ has at most two connected components and they conjectured that the same result holds in every dimension. Moreover, they proved that the result also holds in higher dimensions assuming the termination of klt flips. Later, Birkar (\cite{Bir24}) and Filipazzi-Svaldi (\cite{FS23}) independently proved that if $-(K_X+B)$ is nef, then $\Nklt(X,B)$ has at most two connected components. In fact, the same result holds for generalized pairs. \\
In the last few years, minimal model program (MMP) has been developed for normal compact K\"ahler varieties. In this article, following Birkar's approach, we investigate the corresponding question for proper morphisms between normal compact K\"ahler varieties. The result is already known in a very special situation. Let $(X,B)$ be a compact K\"ahler dlt pair of dimension at most $3$ over a projective variety $S$ satisfying $K_X+B \sim_{\mbQ, f} 0 $. Then $\Nklt(X,B)\cap f^{-1}(s)=\lfloor B \rfloor \cap f^{-1}(s)$ has at most two connected components (see \cite[Proposition 4.1]{SD26}). In this article, we prove the statement in every dimension for generalized pairs whose nef part is allowed to be transcendental. Our main theorem is the following:
\begin{theorem}\label{thm:main}
    Let $(X,B+\beta)$ be a generalized pair over $S$ where $f: X \to S$ be a morphism of normal compact K\"ahler varieties with connected fibers. Assume that $-(K_X+B+\beta)$ is nef over $S$. Then the following holds:
    \begin{enumerate}
        \item If $-(K_X+B+\beta)$ is big over $S$, then $\Nklt(X,B+\beta)\cap f^{-1}(s)$ is connected for every $s\in S$.
        \item If $\Nklt(X,B+\beta)\cap f^{-1}(s)$ is disconnected for some $s\in S$, then $\Nklt(X,B+\beta)\cap f^{-1}(s)$ has exactly two connected components. Moreover, the pair $(X,B+\beta)$ is generalized log canonical over $s$ and there exists a higher model $W$ of $X$ with $K_W+B_W+\beta_W=\phi^*(K_X+B+\beta)$ and there is a contraction $\theta: W \to Z$ over $S$ where $\dim Z=\dim X-1$ such that for any general fiber $F$ of $W\to Z$, $\beta_W|_F\equiv 0$ and $(F,B_W|_F)\cong (\mathbb{P}^1, p_1+p_2)$.
    \end{enumerate}
\end{theorem}
In recent years, significant progress has been done in geometry of projective morphisms between complex analytic spaces (possibly non compact). So it is also very natural to ask whether the similar result holds in this setting. In \cite[Lemma 4.1]{Fu26}, Fujino proved that if $-(K_X+\Delta)$ is $f$-nef and $f$-big where $f:X \to S$ be a projective morphism between complex analytic spaces (possibly non compact) with connected fibers, then $\Nklt(X,\Delta)\cap f^{-1}(s)$ is connected for every $s\in S$. Later the same result has been extended to generalized pair setting by Hashizume (see \cite[Lemma 2.16]{Has26}). We investigate the question assuming only $-(K_X+\Delta)$  is relatively nef.  We record the corresponding connectedness statement for projective morphisms between analytic varieties (see Theorem \ref{thm:con-proj}).\\ 

\textbf{Comparison with Birkar's proof:} Our main Theorem \ref{thm:main} follows from Theorems \ref{thm:main1} and \ref{thm:main2}. We are not able to follow Birkar's argument directly. His argument is local over the base. Given a point of the base over which the non-klt locus is disconnected, he replaces the base by a neighborhood (in fact, an \'etale neighborhood) of this point. In this way he may assume that the non-klt locus itself is disconnected, and then he runs MMP over this small base. In our setting, $f:X \to S$ is only a proper morphism between compact K\"ahler varieties, and it need not be projective. The MMP is only available for such morphism when the base is compact \cite{HX26}. Therefore we cannot replace by a small neighborhood of $s$, and every MMP is run over the whole base $S$. This has two consequences. First, $\Nklt(X,B+\beta)$ may be connected even though $\Nklt(X,B+\beta)\cap f^{-1}(s)$ is not, so we can only use the disconnectedness of one fiber. Second, an MMP over $S$ also changes $X$ away from $f^{-1}(s)$, and it may change the fiber over $s$ itself, so we need to control this fiber through every step of the MMP. On the other hand, for projective morphisms between complex analytic spaces the MMP is available over a neighborhood of a suitable compact subset of the base \cite{fuj}; this is why in Section $4$ we can follow Birkar's local argument more closely. Moreover, in \cite{Bir24}, the nef part of a generalized pair is a b-$\mbR$-divisor. However, in our setting $\beta$ is a b-$(1,1)$-current and its class is in  general not represented by any $\mbR$-divisor. Hence the arguments that need a divisor in the class of a nef part does not work here. So we use different arguments to achieve our desired result. Also, every projective variety is a compact K\"ahler variety. Therefore the Theorem \ref{thm:main} applies to a morphism with connected fibers between normal projective varieties and to a generalized pair $(X,B+\beta)$ whose nef part $\beta$ is b-$(1,1)$-current. The class of $\beta$ need not be the class of an $\mbR$-divisor. Such generalized pairs are not covered by \cite{Bir24}. Thus, to the best of our knowledge, Theorem \ref{thm:main} is new even for projective generalized pairs when the nef part is transcendental.\\

\textbf{AI Disclosure:} Opus 5.5 was used to check the article for typographical errors and to improve the readability of some sentences. All mathematical ideas, arguments, proofs and results in this article are the author's own. The author takes full responsibility for the content of the article.\\

\textbf{Acknowledgements:} The author would like to thank his advisor Prof. Omprokash Das for very useful discussions.

\section{Preliminaries}
We follow the usual conventions of the Minimal Model Program. We refer the reader to the preliminaries of \cite{HX26} for standard notation, definitions and terminology. For the convenience of the reader, we recall the definitions of pair and generalized pairs. 
\begin{definition}
    A pair consists of a normal analytic variety $X$, together with an effective $\mbR$- divisor $B$ such that $K_X+B$ is $\mbR$-Cartier. Let $\phi: Y \to X$ be a log resolution of the pair $(X,B)$. For any prime divisor $D$ on $Y$, we define log discrepancy of $D$ with respect to the pair $(X,B)$ as $a(D, X,B):=- \mult_D(B_Y)$ where $K_Y+B_Y=\phi^*(K_X+B)$. We call $(X,B)$ a \textit {klt} pair (resp. an \textit {lc} pair) if $a(D,X,B)> -1$ (resp. $a(D,X,B)\geq -1$) for every prime divisor $D$ on $Y$. A \textit{non-klt place} (resp. \textit{non-lc place})  is a prime divisor $D$ on a bimeromorphic model $X'$ of $X$ such that $a(D,X, B)\leq -1$ (resp. $a(D,X,B)<-1$). A \textit{non-klt center} (resp. \textit{non-lc center}) is the image of some non-klt place (resp. non-lc place). 
\end{definition}
\begin{definition}
Let $X$ be a normal compact analytic variety and let $f:X \to S$  be a proper morphism where $S$ is relatively compact. Let $\beta$-be a b-(1,1)-current. Let $B\geq 0$ be an $\mbR$-divisor on $X$. We say that $(X,B+\beta)$ be a generalized pair over $S$ if the followings hold:
\begin{itemize}
    \item there is a log resolution $\phi: X' \to X$  such that $\beta_{X'}$ descends on $X'$.
    \item $[\beta_{X'}]$ is nef over $S$.
    \item $[K_{X'}+B_{X'}+\beta_{X'}]=\phi^*[K_X+B+\beta]$ where $\phi_*B_{X'}=B$.
\end{itemize}

We say that $(X,B+\beta)$ is \textit{generalized klt} (resp. \textit{generalized lc}) if the generalized discrepancy $a(D,X,B+\beta)=-\mult_D(B_{X'})>-1$ (resp. $a(D,X,B+\beta)=-\mult_D(B_{X'})\geq -1$) for all prime divisors $D$ over $X$. A \textit{ generalized non-klt place} (resp. \textit{ generalized non-lc place}) is a divisor $D$ on a bimeromorphic model $Y$ such that $a(D,X, B+\beta)\leq -1$ (resp. $a(D,X,B+\beta)<-1$. A \textit{non-klt center} (resp. \textit{non-lc center}) is image of some non-klt place (resp. non-lc place). The \textit{generalized non-klt} $\Nklt (X,B+\beta)$ (resp.\textit{ generalized non-lc} $\Nlc (X,B+\beta)$) locus is union of all \textit{ generalized non-klt centers} (resp. \textit{ generalized non-lc centers}).
\end{definition}
In this article, we write $\alpha\equiv \beta$ if they represent the same class in $H^{1,1}_{BC}(X)$. Also, we say $\alpha\equiv_{S} \beta$ if $\alpha-\beta=f^*\gamma$ where $[\gamma]\in H^{1,1}_{BC}(S)$.
The following lemmas are known to experts; for the reader's convenience, we include their proofs.
\begin{lemma}\label{lem:conn}
Let $X$ be a normal projective variety over $\mbC$ of dimension at least $2$, and let $H$ be an effective ample integral divisor on $X$. Then $\Supp H$ is connected.
\end{lemma}
\begin{proof}
For a positive integer $m$,
consider the following exact sequence of $\mcO_X$-modules,
\[0 \to \mcO_X(-mH)\to \mcO_X \to \mcO_{mH}\to 0\]
From the above exact sequence we obtain a long exact sequence in cohomology,
\[0\to H^0(X,\mcO_X(-mH))\to H^0(X,\mcO_X)\to H^0(X,\mcO_{mH})\to H^1(X, \mcO_X(-mH))\to\cdots\]
Note that $\Supp H= \Supp (mH)$.
By Serre duality and Serre vanishing $H^1(X, \mathcal{O}_X(-mH))=0$ for $m>>0$, and hence we obtain an exact sequence
\[ H^0( \mcO_X)\to H^0(\mcO_{mH})\to 0\].
This implies that $H^0(\mcO_{mH})=\mbC$ and hence $\Supp H$ is connected.
\end{proof}
\begin{lemma}\label{lem:dominating}
Let $f:W\to W'$ be a proper surjective morphism between analytic varieties with $W'$ normal. Let $V\subset W'$ be any connected open subset. Then $f(Q)=V$ for every non empty clopen subset $Q \subset f^{-1}(V)$.
\end{lemma}
\begin{proof}
    Since $Q$ is a clopen subset of $f^{-1}(V)$, it is an analytic subset. As $f|_Q:Q \to V$ is also proper, $f(Q)$ is a closed analytic subset of $V$. Now since $W$ is irreducible, we can find a dense Zariski open subset $V'$ such that $f|_{V'}$ is a submersion. Since $Q$ is open, $V' \cap Q \neq \emptyset$. Thus $f(Q)$ contains a non empty open subset of $V$. Note that $V$ is irreducible as it is connected and normal. Therefore $f(Q)=V$.
\end{proof}
\section{Proof of Theorem \ref{thm:main}}
In this section, we prove our main theorem \ref{thm:main}. 
\begin{theorem}\label{thm:main1}
    Let $(X,B+\beta)$ be a generalized pair over $S$ where $f: X \to S$ is a proper surjective morphism with connected fibers between normal compact K\"ahler varieties. Assume that $-(K_X+B+\beta)$ is nef over $S$. Then the fibers of $\Nklt(X,B+\beta) \to S$ are connected if one of the following conditions holds:
    \begin{enumerate}
        \item $-(K_X+B+\beta)$ is big over $S$.
        \item The restricted morphism $\Nklt(X,B+\beta) \to S$ is not surjective.
    \end{enumerate}
\end{theorem}
\begin{proof}
    \begin{enumerate}
        \item \label{kol-shok-ver}
        By \cite[Theorem 3.3]{HLX26}, we can find a dlt modification of the given generalized pair. That is, we obtain a projective bimeromorphic morphism $\mu:Y \to X$ such that $Y$ is strongly $\mbQ$- factorial, $(Y,\Theta_Y)$ is a dlt pair, satisfying $K_Y+B_Y+\beta_Y=\mu^*(K_X+B+\beta)$ and $B_Y=\Theta_Y+G_Y$ with $\Supp G_Y \subset \lfloor\Theta_Y\rfloor$. Observe that $-(K_Y+B_Y+\beta_Y)$ is nef and big over $S$. Since $-(K_Y+B_Y)$ is relatively big over $S$, $f\circ \mu$ is a Moishezon morphism. Hence $f\circ \mu$ is a projective morphism by \cite[Theorem 3.1]{CH24}. Let $\mathscr{I}$ be the generalized multiplier ideal of $(Y,B_Y+\beta_Y)$. Let $N$ be the closed analytic subset defined by $\mathscr{I}$. By relative Nadel vanishing ( see \cite[Theorem 2.3]{HX26}), $R^1(f\circ\mu)_*\mathscr{I}=0$ for all $i>0$.\\
Applying $(f\circ \mu)_*$ to the short exact sequence $0\to \mathscr{I}\to \mathcal{O}_Y\to \mathcal{O}_Z \to 0$, we obtain a long exact sequence, and hence the map $\mathcal{O}_S=(f\circ \mu)_*\mathcal{O}_Y\to (f\circ \mu)_*\mathcal{O}_Z$ is surjective. By Stein factorization, $Z\cap (f\circ \mu)^{-1}(s)$ is connected for every $s$ for which it is non empty. Since $\Nklt(X,B+\beta)=\mu(Z)$, it follows that $\Nklt(X,B+\beta)\cap f^{-1}(s)=\mu(Z \cap (f\circ \mu)^{-1}(s))$ is connected for every $s\in S$. \\
\item We prove the result in several steps.\\

\textbf{Step 1:} Set $\alpha:=-(K_X+B+\beta)$, which is nef over $S$ by assumption. Then $(X,B+\beta+\alpha)$ is a generalized pair over $S$ and $\alpha$ descends to $X$. Hence, for every prime divisor $D$ over $X$, the generalized discrepancy satisfy $a(D, X, B+\beta)=a(D,X,B+\beta+\alpha)$. 
Replacing $(X,B+\beta)$ by $(X,B+\beta+\alpha)$, we may assume that $K_X+B+\beta\equiv_S 0$.\\

\textbf{Step 2:} Let $\phi: W \to X$ be a log resolution of the pair $(X,B+\beta)$ such that $\beta_W$ descends on $W$. Write $B=\sum c_i C_i+\sum d_j D_j$ where $C_i$ and $D_j$ are distinct prime divisors, $c_i<1$ and $d_j\geq 1$. Set $B^*=\sum c_i C_i +\sum D_j$ and $\Theta_W=\phi_*^{-1}B^*+\sum E_i$ where the $E_i$ are $\phi$-exceptional prime divisors. For any $t\in (0,1)$, $(W, \Theta_W-t\lfloor\Theta_W\rfloor+\beta_W)$ is a generalized klt pair. By \cite[Theorem 1.5]{HX26}, we may run a $(K_W+\Theta_W-t\lfloor\Theta_W\rfloor+\beta_W)$-MMP over $X$ and obtain a morphism $\phi': W' \to X$ satisfying $K_{W'}+\Theta_{W'}-t\lfloor\Theta_{W'}\rfloor+\beta_{W'})$ is nef over $X$. Let $B_{W'}$ be the $\mbR$-divisor defined by $K_{W'}+B_{W'}+\beta_{W'}=\phi'^*(K_X+B+\beta)$, and set $D_t:= B_{W'}-(\Theta_{W'}-t\lfloor \Theta_{W'}\rfloor)$. Then
\[ K_{W'}+B_{W'}+\beta_{W'}=\phi'^*(K_X+B+\beta)\]
\[-D_t=(\Theta_{W'}-t\lfloor\Theta_{W'}\rfloor)-B_{W'}\]
\[=(K_{W'}+\Theta_{W'}-t\lfloor \Theta_{W'}\rfloor+\beta_{W'})-(K_{W'}+B_{W'}+\beta_{W'})\]
\[\equiv_X (K_{W'}+\Theta_{W'}-t\lfloor \Theta_{W'}\rfloor+\beta_{W'}).\]
Thus $-D_t$ is $\phi'$-nef. Note that $\phi'_*D_t$ is effective. Thus by the negativity lemma (see \cite[Lemma 1.3]{Wang21}), $D_t$ is effective. For a prime divisor $E$ on $W$, let $b_E$ denote the coefficient of $E$ in $B_W$, where $K_W+B_W+\beta_W=\phi^*(K_X+B+\beta)$.\\

\begin{claim}\label{clm:b_E1} We can choose $t\in (0,1)$ sufficiently small such that on the final model $W'$ the following holds:
    if $E$ is a $\phi$-exceptional divisor on $W$ which survives on $W'$, then $b_E\geq 1$.
\end{claim}
\begin{proof}[Proof of Claim \ref{clm:b_E1}]
    We may choose sufficiently small $t$ such that $0<t<1-b_E$ for every exceptional divisor $E$ with $b_E<1$. Run a $(K_W+\Theta_{W}-t\lfloor \Theta_{W}+\beta_{W})$-MMP over $X$, and let $\phi':W \to X$, be the  resulting minimal model over $X$. Let $E$ be a $\phi$-exceptional divisor on $W$ with $b_E<1$ which survives on $W'$. Then $\mult_E(D_t)=b_E-(1-t)<0$ which contradicts the fact that $D_t$ is effective. Thus $b_E\geq 1$ for every exceptional divisor $E$ which survives on $W'$.
\end{proof}

Set $G_{W'}:=B_{W'}-\Theta_{W'}$. For any $\phi$-exceptional divisor $E$ on $W'$, we have $\mult_E(G_{W'})=b_E-1\geq 0$, and for strict transform of any component of $B$, we have $\mult_E(G_{W'})=\mult_E(B_W)-\mult_E(\Theta_W)\geq 0$. Thus $G_{W'}$ is effective and $\Supp G_{W'}\subseteq \lfloor \Theta_{W'}\rfloor$. Replacing $(X,B+\beta)$ by $(W',B_{W'}+\beta_{W'})$, we may assume that the following condition holds:

$(*) B=\Theta+G$ where $G \subseteq \Supp \lfloor\Theta\rfloor)$ such that $(X, \Theta-t\lfloor\Theta\rfloor+\beta)$ is a strongly $\mbQ$- factorial generalized klt pair for some $t\in (0,1)$.\\

This implies that $\Supp\lfloor\Theta\rfloor \subseteq \Nklt(X,B+\beta)\subseteq \Supp (t\lfloor\Theta\rfloor+G)= \Supp\lfloor\Theta\rfloor$ and hence $\Nklt(X,B+\beta)=\lfloor\Theta\rfloor$. Since $\Nklt(X,B+\beta)\to S$ is not surjective, $\lfloor\Theta\rfloor$ is vertical over $S$. Thus we have $K_X+\Theta-t\lfloor\Theta\rfloor+\beta \equiv_U 0$ for some nonempty Zariski open subset $U\subset S$. Hence $K_X+\Theta-t\lfloor\Theta\rfloor+\beta$ is pseudo effective over $S$.\\

\textbf{Step 3:} We can run a $(K_X+\Theta-t\lfloor\Theta\rfloor+\beta)$-MMP over $S$ with scaling of a K\"ahler class $\omega$. The MMP need not terminate. However, if $\lambda_i$ denotes the scaling number at the $i$-th step $X_i$ of this MMP, then $\lambda_i \to 0$. This follows from finiteness of weak log canonical models (see \cite[Theorem 5.15]{HX26}). Note that since $K_X+\Theta+G+\beta\equiv_S 0$, the condition  $(*)$ also holds at every step of the MMP. Moreover, since $K_X+\Theta-t\lfloor\Theta\rfloor+\beta$
is pseudo effective over $S$, every step of the MMP is either a flip or a divisorial contraction. \\

\textbf{Step 4:} In this step, we show that the number of connected components of $\Nklt(X,B+\beta) \cap f^{-1}(s)$ remains the same at every step of the MMP. Let $u: X \to X'$ be a divisorial contraction over $S$, let $g:X' \to S$ be the induced morphism, and let $B'$ and $\beta'$ denote pushforwards of $B$ and $\beta$ to $X'$. Note that $K_X+B+\beta\equiv_{X'}0$. Since $X \to X'$ is bimeromorphic, $-(K_X+B+\beta)$ is nef and big over $X'$. Then by \ref{kol-shok-ver}, $\Supp (t\lfloor\Theta\rfloor+G)\cap u^{-1}(x')$ is connected for every pt $x'\in X'$. Since $u(\Nklt(X,B+\beta))=\Nklt(X',B'+\beta')$, the sets $\Nklt(X,B+\beta)\cap f^{-1}(s)$ and $\Nklt(X',B'+\beta')\cap g^{-1}(s)$ have the same number of connected components. Let $\psi:X\dashrightarrow X'$ be a flip, with flipping contraction $u:X \to V$ and flipped contraction $u^+:X'\to V$.
Consider the following commutative diagram:

\[\begin{tikzcd}
	X && {X'} \\
	& V \\
	& S
	\arrow["\psi", dashed, from=1-1, to=1-3]
	\arrow["u", from=1-1, to=2-2]
	\arrow["f"', from=1-1, to=3-2]
	\arrow["{u^+}"', from=1-3, to=2-2]
	\arrow["g", from=1-3, to=3-2]
	\arrow["h"', from=2-2, to=3-2]
\end{tikzcd}\]
By a similar argument, $\Nklt(X,B+\beta)\cap u^{-1}(v)$ and $\Nklt (X',B'+\beta')\cap (u^+)^{-1}(v)$ are connected for every point $v\in V$. Hence we obtain the desired result.\\
\textbf{Step 5:} Suppose, for contradiction that the fiber of $\Nklt(X,B+\beta)\to S$ over $s$ is not connected. In this step we will show that after finitely many steps of the above MMP, the irreducible components of the fiber over $s$ that are not contained in the $\Nklt(X,B+\beta)$ stabilize. Let $F_1,F_2,...,F_k$ be the irreducible components of $f^{-1}(s)$. Since $f$ has connected fibers and the fiber of $\Nklt(X,B+\beta)\to S$ over $s$ is not connected, we have $f^{-1}(s)\nsubseteq \Nklt(X,B+\beta)$. Thus some $F_i$ is not be contained in $\Nklt(X,B+\beta)$. Reordering the indices, we may assume that $F_1,...,F_m$ are not contained in $\Nklt(X,B+\beta)$ and that the other components are contained in it. \\
Let $\psi:X\dashrightarrow X'$ be the first step of the MMP, and suppose that $\psi$ is a divisorial contraction.\\
\begin{claim}\label{clm:im-exp-in-Nklt}
Let $E$ be the exceptional divisor of $\psi$. Then $\psi(E)\subset \Nklt(X',B'+\beta')$.
\end{claim}
\begin{proof}[Proof of the Claim \ref{clm:im-exp-in-Nklt}]
    If not, choose a point $p\in \psi(E)\setminus \Nklt(X',B'+\beta')$. Since $\psi$ is a projective morphism, we can choose a curve $\Gamma\subset \psi^{-1}(p)$. Then $(t\lfloor\Theta\rfloor+G)\cdot \Gamma=0$ which contradicts the fact that $(t\lfloor\Theta\rfloor+G)$ is $\psi$-ample. 
\end{proof}
Similarly, we can show that if $\psi: X \dashrightarrow X'$ is a flip, then the flipped locus is also contained in $\Nklt(X',B'+\beta')$. This shows that any irreducible component of the fiber $X'\to S$ over $s$ that is not contained in $\Nklt(X',B'+\beta')$ is the bimeromorphic pushforward of some $F_i$ with $1\leq i \leq m$. Thus, after finitely many steps, replacing $X$ by the model appearing at that step and possibly decreasing $m$, we may assume that, on every subsequent model of the MMP, the bimeromorphic pushforwards of $F_1, F_2,...,F_m$ are exactly the irreducible components of the fiber over $s$ not contained in $\Nklt(X',B'+\beta')$.\\

\textbf{Step 6:}
Since $\Nklt(X, B+\beta)\cap f^{-1}(s)$ is disconnected, there exists an irreducible component $F_i (1\leq i\leq m)$ such that $F_i \cap \Nklt(X,B+\beta)\neq \emptyset$. We denote this component by $F$. Let $X:=X_0 \dashrightarrow X_1 \dashrightarrow X_2\dashrightarrow \cdots$ be the steps of the given MMP, let $\phi_j : X \dashrightarrow X_j$ be the induced  bimeromorphic map, and let $B_j,\Theta_j,\omega_j,\beta_j$ denote the pushforwards of $B,\Theta,\beta$ and $\omega$ to $X_j$. Since $(\phi_j)_* F\nsubseteq \Nklt(X_j, B_j+\beta_j)$ and $\Nklt(X_j, B_j+\beta_j)$ contains the image of the $\phi_j$ exceptional locus, $\phi_j$ is an isomorphism near a general point of $F$. Let $W$ be a common log resolution which resolves the map $\phi_j$. Let $\bar{F}$ be the strict transform of $F$ on $W$. Since $\phi_j$ is $(K_X+\Theta-t\lfloor\Theta\rfloor+\beta_X+\lambda_j \omega)$-non positive and $\phi_j^{-1}$ has no exceptional divisor, we can write:
\[p^*(K_X+\Theta-t\lfloor\Theta\rfloor+\beta_X+\lambda_j \omega)=q^*(K_{X_j}+\Theta_j-t\lfloor\Theta_j\rfloor+\beta_{X_j}+\lambda_j\omega_j)+ E\]
where $E$ is effective $q$-exceptional divisor. Note that $\bar{F} \nsubseteq \Supp E$. Thus $E|_{\bar{F}}$ is an effective divisor on $\bar{F}$. Restricting the above equation to $\bar{F}$, we obtain,
\[p^*((K_X+\Theta-t\lfloor\Theta\rfloor+\beta_X+\lambda_j \omega)|_F)=q^*((K_{X_j}+\Theta_j-t\lfloor\Theta_j\rfloor+\beta_{X_j}+\lambda_j\omega_j)|_{(\phi_j)_*F})+ E|_{\bar{F}}\]
As $(K_{X_j}+\Theta_j-t\lfloor\Theta_j\rfloor+\beta_{X_j}+\lambda_j\omega_j)$ is nef over $S$, $(K_{X_j}+\Theta_j-t\lfloor\Theta_j\rfloor+\beta_{X_j}+\lambda_j\omega_j)|_{(\phi_j)_*F}$ is nef. Thus right hand side is pseudo effective on $\bar{F}$
This shows that $(K_X+\Theta-t\lfloor\Theta\rfloor+\beta_X+\lambda_j \omega)|_F$ is pseudo effective for every scaling number $\lambda_j$ appearing in the above MMP.\\

\textbf{Step 7:}
Set $D_F=(t\lfloor\Theta\rfloor+G)|_F$. Note that $D_F$ is a non zero effective divisor on $F$, since $F\cap \Nklt(X,B+\beta)\neq \emptyset$ and $F\nsubseteq \Nklt(X,B+\beta)$. Now we have $((K_X+\Theta-t\lfloor\Theta\rfloor+\beta_X)|_F\equiv -D_F$, since $(K_X+\Theta+G+\beta_X)\equiv_S0$. 
Fix a K\"ahler class $\omega_F$ on $F$. Then $(K_X+\Theta-t\lfloor\Theta\rfloor+\beta_X+\lambda_j \omega)|_F\cdot \omega_F^{\dim F-1}\geq 0$. As $\lambda_j\to 0$, we get $(K_X+\Theta-t\lfloor\Theta\rfloor+\beta_X)|_F\cdot \omega_F^{\dim F-1}\geq 0$. But 
$(K_X+\Theta-t\lfloor\Theta\rfloor+\beta_X)|_F\cdot \omega_F^{\dim F-1}=-D_F\cdot\omega_F^{\dim F-1}<0$ which gives a contradiction. This proves that the fibers of $\Nklt(X,B+\beta)\to S$ are connected.

\end{enumerate}

\end{proof}

Now we consider the case when fiber of $\Nklt(X,B+\beta)\to S$ over $s\in S$ is not connected.
\begin{theorem}\label{thm:main2}
    Let $(X, B+\beta)$ be a generalized pair over $S$ where $f:X \to S$ is a proper surjective morphism with connected fibers between normal compact K\"ahler varieties. Assume that $-(K_X+B+\beta)$ is nef over $S$ and that the fiber of $\eta: \Nklt(X,B+\beta) \to S$ is not connected for some $s\in S$. Then the followings hold:
    \begin{enumerate}
\item  the restricted morphism $\eta: \Nklt(X,B+\beta)\to S$ is surjective and its fiber over $s$ has exactly two connected components.
  \item the pair $(X,B+\beta)$ is generalized log canonical over $s$, and there exist a higher model $W$ of $X$ with $K_W+B_W+\beta_W=\phi^*(K_X+B+\beta)$ and there is a contraction $\theta: W \to Z$ over $S$ where $\dim Z=\dim X-1$ such that for any general fiber $F$ of $W\to Z$, $\beta_W|_F\equiv 0$ and $(F,B_W|_F)\cong (\mathbb{P}^1, p_1+p_2)$.
 \end{enumerate}
\end{theorem}

\begin{proof}
We prove $(1)$ in several steps. \\

\textbf{Step 1:}
As in Step 1 of the proof of Theorem \ref{thm:main1} Step 1, we may assume that $K_X+B+\beta \equiv_S0$. By Step 2 of the proof of Theorem \ref{thm:main1}, we can replace  $(X,B+\beta)$ by a generalized pair $(Y,B_Y,\beta_Y)$ and assume that the condition $(*)$ holds:\\
 $(*) B=\Theta+G$ where $G \subseteq \Supp \lfloor\Theta\rfloor)$ such that $(X, \Theta-t\lfloor\Theta\rfloor+\beta)$ is a strongly $\mbQ$- factorial generalized klt pair for some $t\in (0,1)$. In particular, we have $\Nklt(X,B+\beta)= \Supp (t\lfloor\Theta\rfloor+G)=\lfloor\Theta\rfloor$.\\
 
 \textbf{Step 2:}
 Since there exists $s\in S$ such that fiber of $\eta$ is disconnected, $\eta$ is surjective by Theorem \ref{thm:main1}(2). Note that $K_X+\Theta-t\lfloor\Theta\rfloor+\beta\equiv_S -(t\lfloor\Theta\rfloor+G)$ is not pseudo effective over $S$, since $t\lfloor\Theta\rfloor+G)$ dominates $S$. We may run a $(K_X+\Theta-t\lfloor\Theta\rfloor+\beta)$-MMP over $S$ and obtain a Mori fiber space $u: X' \to Z$ over $S$. We obtain the following commutative diagram:

\[\begin{tikzcd}
	X && {X'} & \\
	&&& Z \\
	& S
	\arrow["\psi", dashed, from=1-1, to=1-3]
	\arrow["f"', from=1-1, to=3-2]
	\arrow["u", from=1-3, to=2-4]
	\arrow["g"', from=1-3, to=3-2]
	\arrow["h", from=2-4, to=3-2]
\end{tikzcd}\]
Let $B',\Theta',G',\beta'$ denote the push forwards of $B,\Theta, G$ and $\beta$ to $X'$. Note that $t\lfloor\Theta'\rfloor+G'$ is ample over $Z$. Hence there exists an irreducible component $\Lambda\subset \Nklt(X',B'+\beta')$ such that $S$ is $u$-ample. \\

\textbf{Step 3:} We show that $u$ has relative dimension one.
\begin{claim}\label{clm: rel-dim-1}
    $\dim X'-\dim Z=1$.
\end{claim}
\begin{proof}[Proof of the Claim \ref{clm: rel-dim-1}]
If possible, assume that $\dim X'-\dim Z=r\geq 2$. For a general fiber $F$ of $u$, $\Lambda|_F$ is an ample $\mbQ$- Cartier divisor on $F$. Thus $S|_F$ is connected by Lemma \ref{lem:conn}. By Stein factorization we obtain the following:
\[\begin{tikzcd}
	\Lambda && {Z_{\Lambda}} \\
	& Z
	\arrow["{\eta_1}", from=1-1, to=1-3]
	\arrow["{u|_S}"', from=1-1, to=2-2]
	\arrow["{\eta_2}", from=1-3, to=2-2]
\end{tikzcd}\]
\end{proof}
Here $\eta_1$ has connected fibers and $\eta_2$ is a finite map. Now, since $\Lambda|_{u^{-1}(z)}$ is connected for a general point $z\in Z$, $\eta_2$ is bimeromorphic. Since $Z$ is normal, $\eta_2$ is an isomorphism by Zariski's main theorem. Therefore $\Lambda|_{u^{-1}(z)}$ is connected for every point $z\in Z$. If $\Lambda'$ is any other irreducible component of $\Nklt(X',B'+\beta')$ that dominates $Z$, then $\Lambda'$ is also $u$-ample, and hence $\Lambda'|_{u^{-1}(z)}$ is also connected for every point $z\in Z$. Since $h$ has connected fibers, $h^{-1}(s)$ is connected. This implies that $\Lambda_s:=\Lambda\cap g^{-1}(s)$ is connected for every $s\in S$. Now, for any other horizontal component $\Lambda'$, the set $\Lambda'\cap u^{-1}(z)$ is positive dimensional. Thus $\Lambda_s\cap \Lambda'_s\neq \emptyset$. Let $T$ be any component of $\Nklt(X',B'+\beta')$ that is vertical with respect to $u$. Then every connected component of $T_s:=T \cap g^{-1}(s)$ intersects $\Lambda_s$. This shows that $\Nklt(X',B'+\beta')\cap g^{-1}(s)$ is connected, which implies by Step (4) of the proof of Theorem \ref{thm:main1}, that $\Nklt(X,B+\beta)\cap f^{-1}(s)$ is also connected. This gives a contradiction. 
\end{proof}
Since $u$ is a Mori fiber space contraction of relative dimension $1$, a general fiber $F$ of $u$ is isomorphic to $\cong \mathbb{P}^1$.\\

\textbf{Step 4:} 
Since $u$ is $(K_{X'}+B'+\beta')$-trivial, $(K_{X'}+B'+\beta')|_F\equiv 0$. Let $S_1,...,S_p$ be the irreducible components of $\lfloor \Theta'\rfloor$ that are horizontal with respect to $u$. Let $R_1,....,R_q$ be the other irreducible components of $B'$ that are horizontal with respect to $u$.
Let $d_i$ and $e_j$ be the degree of $u|_{S_i}: S_i \to Z$ and $u|_{R_j}:R_j\to Z$ respectively. Let $a_i$ and $b_j$ denote the coefficients of $S_i$ and $R_j$ respectively in $B'$. Note that $a_i\geq 1$ and $b_j\in (0,1)$ for all $1\leq i \leq p$ and $1\leq j\leq q$. Then we have $\deg (K_F+\beta'_F)+\sum d_ia_i+\sum e_jb_j=0$.
Now $\deg(\beta'_F)\geq 0$, since $\beta'$ is pseudo effective over $S$. This implies that $\sum d_ia_i+\sum e_jb_j \leq 2$. Note that there exists at least one horizontal component $\Lambda \subset \lfloor\Theta'\rfloor$ where $B'=\Theta'+G'$.
There are following three possibilities:
\begin{enumerate}
    \item $p=1$ and $ d_1=1$. Then $a_1\in [1,2] $ and $\deg(\beta'_F) +\sum e_jb_j=2-a_1$.
In this case $S_1$ is the only horizontal component of $\Nklt(X',B'+\beta')$ and all the other components are vertical. Then by same argument as above we can say that $\Nklt(X',B'+\beta')\cap g^{-1}(s)$ is connected for every $s\in S$ which is a contradiction.
\item $p=2$, $d_1=d_2=1$. Then $a_1=a_2=1, q=0$ and $\deg\beta'_F=0$.\label{itm:2}
In this case we obtain $(F,B'+\beta'_F)\cong (\mathbb{P}^1, p_1+p_2)$ and $\beta'_F\equiv 0$.
\item $p=1, d_1=2$. Then $a_1=1$ and $q=0, \deg\beta'_F=0$\label{itm:3}
In this case, we also obtain $(F,B'_F)\cong (\mathbb{P}^1, p_1+p_2)$ and $\beta'_F\equiv 0$.\\
\end{enumerate}

\textbf{Step 5:} In this step, we show that number of connected components of $\Nklt(X,B+\beta)\cap f^{-1}(s)$ is at most $2$.
It is enough to show that the number of connected components of $\Nklt(X',B'+\beta')\cap g^{-1}(s)$ is at most $2$. 
In case \ref{itm:2}, there is no component $T$ of $\Nklt(X',B'+\beta')$ that is vertical with respect to $u$ and meets $g^{-1}(s)$. Indeed, otherwise every connected component of $T_s$ would intersect both $\Lambda_s$ and $\Lambda'_s$, and hence $\Nklt(X',B'+\beta')\cap g^{-1}(s)$ would be connected, which is a contradiction. Thus $\Nklt(X',B'+\beta')\cap g^{-1}(s)= \Lambda_s \sqcup \Lambda'_s$ whenever $\Nklt(X',B'+\beta')\cap g^{-1}(s)$ is disconnected.
In case \ref{itm:3}, there is only one horizontal component $\Lambda$, and it has degree $2$ over $Z$. If $T$ is any component of $\Nklt(X',B'+\beta')$ that is vertical with respect to $u$ and meets $g^{-1}(s)$, then every connected component of $T_s$ meets $\Lambda_s$. Thus $\Nklt(X',B'+\beta')\cap g^{-1}(s)$ and $\Lambda_s$ have the same number of connected components. Let $\bar{\Lambda}\to \Lambda$ be the normalization of $\Lambda$. By Stein factorization, we obtain the following commutative triangle:
\[\begin{tikzcd}
	{\bar{\Lambda}} && {Z'} \\
	& Z
	\arrow[from=1-1, to=1-3]
	\arrow["\bar{u}", from=1-1, to=2-2]
	\arrow[from=1-3, to=2-2]
\end{tikzcd}\]
Now $Z'\to Z$ is a finite map of degree two. Therefore there exists an involution $\sigma:Z' \to Z'$ such that $Z\cong Z'/\langle\sigma\rangle$. Note that $h^{-1}(s)$ is connected, and hence $\bar{u}^{-1}(h^{-1}(s))$ has at most two connected components. Therefore $S_s$ has at most two connected components [and by the standing assumption that the fiber is disconnected, exactly two].\\
Now we prove $(2)$.\\
Since $S$ is normal, we can choose an irreducible open neighborhood $U\subset S$ of $s$, such that $\Nklt(X',B'+\beta')\cap g^{-1}(U)= N_1 \sqcup N_2$, where $N_i\cap g^{-1}(s)=C_i$ and $\Nklt(X',B'+\beta')\cap g^{-1}(s)=C_1 \sqcup C_2$. 
Note that since $h$ has connected fibers and $U$ is connected, $h^{-1}(U)$ is also connected. Since $Z$ is normal, $h^{-1}(U)$ is irreducible. For $N_i$ we can find an horizontal component $S_i$ such that $S_i \cap N_i\neq \emptyset$ ( $S_1$ and $S_2$ are not necessarily distinct). 
\begin{claim}\label{clm:Ndominates}
    $N_i$ dominates $h^{-1}(U)$ for both $i=1,2$.
\end{claim}
\begin{proof}[Proof of the Claim \ref{clm:Ndominates}]
Fix $i\in \{1,2\}$ and set $N:=N_i$. Let $\Lambda$ be the horizontal component such that $S\cap N\neq \emptyset$. Let $u|_\Lambda: \Lambda \to Z$ be the projective surjective morphism. Applying Lemma \ref{lem:dominating}, with $ W:=\Lambda, W':=Z, f:=u|_\Lambda, Q:=\Lambda\cap N, V:=h^{-1}(U)$, we obtain the desired result.
\end{proof}

Write $D=D_1+D_2$ where $ D_i$ supports $N_i$.\\
\begin{claim}\label{clm:D-nef-big}
    $D_i$ is relatively nef and big over $h^{-1}(U)$ with respect to $u$.
\end{claim}
\begin{proof}[Proof of the Claim \ref{clm:D-nef-big}]
 Let $\Gamma$ be a curve contracted by $u$ over $h^{-1}(U)$. If $\Gamma \nsubseteq \supp D_i$, then $D_i\cdot \Gamma\geq 0$. If $\Gamma \subset \Supp D_i$, then $D_j\cdot \Gamma=0$ and hence $D_i\cdot \Gamma=D\cdot \Gamma>0$. Since relative dimension is $1$, and $D_i$ has positive degree on general fiber $F$ of $u$ over $h^{-1}(U)$, each $D_i$ is also relatively big.
\end{proof}
We show that $(X',B'+\beta')$ is generalized log canonical over $s$. Suppose, for contradiction, that this is not the case. Then there is a divisor $E$ over $g^{-1}(U)$ with $a(E,g^{-1}(U), B'+\beta')<-1$ whose center meets $g^{-1}(s)$.  Without loss of generality, we may assume that its center meets $N_1$. Over $U$, set $B'_\epsilon=B'-\epsilon D_1$ where $\epsilon>0$ is chosen sufficiently small so that $B'_\epsilon\geq 0$. We have
\[a(E, g^{-1}(U), B'_\epsilon+\beta')=a(E, g^{-1}(U), B'+\beta')+ \epsilon\mult_E(D_1)\]
Taking $\epsilon$ smaller if necessary, we have $a(E, g^{-1}(U), B'_\epsilon+\beta')<-1$. Choose a point $z\in h^{-1}(s)$ such that the image under $u$ of the center of $E$ contains $z$. For any divisor $E'$ over $g^{-1}(U)$ with center contained in $N_2$, we have $a(E', g^{-1}(U), B'_\epsilon+\beta')= a(E, g^{-1}(U), B'+\beta)$. This shows that $\Nklt(g^{-1}(U), B'_\epsilon+\beta')=N_1 \sqcup N_2$.
Now $u|_{g^{-1}(U)}: g^{-1}(U)\to h^{-1}(U)$ is a projective morphism with connected fibers. Moreover, $-(K_{g^{-1}(U)}+B'_\epsilon+\beta')\equiv_{h^{-1}(U)} \epsilon D_1$ is relatively nef and big by Claim \ref{clm:D-nef-big}.  Thus, by an argument similar to the proof of of Theorem \ref{thm:main1}(1), we obtain that $\Nklt(g^{-1}(U), B'_\epsilon+\beta')\cap u^{-1}(z)$ is connected. But $\Nklt(g^{-1}(U), B'_\epsilon+\beta')\cap u^{-1}(z)=(N_1\cap u^{-1}(z))\sqcup (N_2 \cap u^{-1}(z))$ for $z\in h^{-1}(s)$, and both sets on the right are not empty by Claim \ref{clm:Ndominates}. This is a contradiction. \\
Finally, choose a log resolution $W$ of $X$ which resolves the bimeromorphic map $X \dashrightarrow X'$. Since image of the exceptional locus of $W\to X'$ has codimension at least $2$ in $X'$, it cannot dominate $Z$. Thus a general fiber of $W\to Z$ is isomorphic to corresponding fiber of $X'\to Z$, and hence we obtain our desired result.\\

\begin{proof}[Proof of Theorem \ref{thm:main}]
    The result follows from Theorems \ref{thm:main1} and \ref{thm:main2}.
\end{proof}
\section{Connectedness Principle for Projective morphism between normal analytic varieties}
In this section we discuss the connectedness principle for projective morphisms between normal analytic varieties (possibly non compact). We follow same notations and terminologies defined in the subsection named "Generalized pairs" in \cite[Subsection 2.2]{Has26}. In our proof, we run a relative MMP for projective morphisms between complex analytic spaces; the following remark is therefore essential.
\begin{remark}
 In \cite{fuj}, Fujino proved that all BCHM-type result hold for projective morphisms between complex analytic spaces over a neighborhood of a compact subset $W$ which satisfies the property \textbf{(P)}. Moreover, for any point $s\in S$, we can always find a Stein compact neighborhood $W$ of $s$ which satisfies the property \textbf{(P)} (see \cite[Remark 2.18]{DHY26}).
\end{remark}
\begin{remark}\cite[Proposition 4.6]{Fu26}
    The following theorem is already known in a special situation. Let $f:X \to S$ be a projective morphism with connected fibers between normal analytic varieties (possibly non compact). Let $(X,B)$ be a dlt pair such that $K_X+B\sim_{ \mbQ, f} 0$. Then $\Nklt(X,B)\cap f^{-1}(s)=\lfloor B \rfloor \cap f^{-1}(s)$ has at most two connected components for every $s\in S$.
\end{remark}
\begin{theorem}\label{thm:con-proj}
    Let $(X,B+M)$ be a generalized pair over $S$ where $f:X \to S$ is a projective surjective morphism with connected fibers between normal analytic varieties (possibly non compact). Assume that $-(K_X+B+M)$ is nef over $S$. Then, the following hold:
    \begin{enumerate}
        \item If $\Nklt(X,B+M)\to S$ is not surjective, then $\Nklt(X,B+M)\cap f^{-1}(s)$ is connected for every $s\in S$.
        \item If $\Nklt(X,B+M)\cap f^{-1}(s)$ is not connected for some $s\in S$, then it has exactly two connected components. Moreover, the pair $(X,B+M)$ is generalized log canonical over $s$ and, after base change to an analytic neighborhood of $s$, there exist a higher model $\phi:X'\to X$ with contraction $X'\to Y'$ over $S$ which satisfies the following:
        $K_{X'}+B'+M'=\phi^*(K_X+B+M)$ , $M|_{F'}\equiv 0$ and $(F', B'|_{F'})\cong (\mbP^1, p_1+p_2 )$ where $F'$ is general fiber of $X'\to Y'$.
    \end{enumerate}
\end{theorem}
\begin{proof}
The proof follows Birkar's argument \cite[Theorems 3.1 and 3.5]{Bir24}, using corresponding analytic inputs. Given $s\in S$, we can find a Stein compact neighborhood $W\subset U \subset S$ which satisfies property \textbf{(P)}. By \cite[Theorem 2.15]{Has26}, after shrinking $S$ around $W$, we can find generalized dlt modification. Thereafter, we can run the required MMPs and carry out the arguments over neighborhoods of $W$, shrinking the base as required. Every such neighborhood contains $s$, so the resulting conclusion applies to the original fiber over $s$.
\end{proof}
\bibliographystyle{alpha}
\bibliography{references}
\end{document}